\documentclass[a4paper,11pt,onecolumn]{amsart}
\usepackage{amsfonts}
\title[Weak mixing and uniform rigidity]
{On ergodic transformations that are both weakly mixing and uniformly rigid}

\newtheorem{thm}{Theorem}
\newtheorem{lemma}[thm]{Lemma}
\newtheorem{defin}[thm]{Definition}

\newtheorem{cor}[thm]{Corollary}
\newtheorem{prop}[thm]{Proposition}

\newtheorem{quest}[thm]{Question}

\numberwithin{thm}{section}

\newcommand{\eps}{\epsilon}

\newcommand{\bZ}{\mathbb{Z}}
\newcommand{\bR}{\mathbb{R}}

\newcommand{\bN}{\mathbb{N}}
\newcommand{\bT}{\mathbb{T}}

\newcommand{\tR}{\triangle}

\newcommand{\bQ}{\mathbb{Q}}

\begin{document}
\bibliographystyle{plain} 

\author[J. James]{Jennifer James}
\address[Jennifer James]{ Department of Mathematics, Brandeis University\\ Waltham, MA 02454, USA}
\email{ jjames@brandeis.edu}
\author[T. Koberda]{Thomas Koberda}
\address[Thomas Koberda]{ Department of Mathematics, Harvard University\\ Cambridge, MA 02138, USA  }
\email{ koberda@math.harvard.edu}
\author[K. Lindsey]{Kathryn Lindsey}
\address[Kathryn Lindsey]{ Department of Mathematics, Cornell University\\ Ithaca, NY 14853, USA}
\email{ klindsey@math.cornell.edu}
\author[C.E. Silva]{Cesar E. Silva}
\address[Cesar Silva]{Department of Mathematics\\
     Williams College \\ Williamstown, MA 01267, USA}
\email{csilva@williams.edu}
\author[P. Speh]{Peter Speh}
\address[Peter Speh]{ Department of Mathematics, Massachusetts Institute of Technology\\ Cambridge, MA 02139, USA}
\email{ pspeh@math.mit.edu}

\begin{abstract}
We examine some of the properties of uniformly rigid transformations, and analyze the compatibility
of uniform rigidity and (measurable) weak mixing along with some of their asymptotic convergence properties.  We show that on Cantor space, there does not exist a finite measure-preserving, totally ergodic, uniformly rigid transformation.  We briefly discuss general group actions and show that (measurable) weak mixing and uniform rigidity can coexist in a more general setting.
\end{abstract}
\maketitle
\begin{center}
\today
\end{center}
\section{Introduction and preliminaries}
The notion of uniform rigidity was introduced by Glasner and Maon in \cite{MR1007412} in the context of topological dynamics.  It was introduced as a topological version of rigidity.  In their paper, Glasner and Maon show that uniform rigidity has many properties analogous to those of rigidity.

It is a classical theorem of Halmos that in the weak topology on the group of all finite measure-preserving transformations of a fixed Lebesgue space, the weakly mixing transformations form a dense $G_{\delta}$.  It is well known that rigid finite measure-preserving transformations are generic.  (A proof for the nonsingular and infinite measure-preserving cases is given in \cite{MR2032829}.)  In particular, there exist (measurably) weakly mixing rigid transformations.  Glasner and Maon show that topologically, uniform rigidity is generic in the following sense: taking an irrational rotation on the circle and crossing with the conjugacy class of a certain transformation on another compact metric space, there is a generic set of uniformly rigid, minimal, and topologically weak mixing transformations on the product.  Unlike rigidity, uniform rigidity seems to have neither a measurable nor a functional analytic characterization, so a generic class of measurably weakly mixing uniformly rigid transformations seems unlikely to exist.  In fact, we will show that every weakly mixing transformation has a realization that is not uniformly rigid, and that uniform rigidity and weak mixing are mutually exclusive notions on a Cantor set, obtaining:
\begin{thm}\label{t:mrah}
Every totally ergodic finite measure-preserving transformation on a Lebesgue space has a representation that is not uniformly rigid, except in the case where the space consists of a single atom.
\end{thm}
The proof of this theorem will be given at the end of section \ref{s:mrah}.

All spaces will be considered simultaneously as topological spaces and as measure spaces.  All results in this paper will concern either the measurable dynamics on $X$ or the interplay between the measurable and topological dynamics.  Unless otherwise noted, we assume all measures to be positive, Borel, regular, complete and finite on compact sets.  When they are in fact finite, they will be normalized to have total weight $1$.  All transformations are assumed to be measure-preserving and invertible unless otherwise noted.  $C(X)$ denotes the set of continuous functions on a space $X$, and $L^p(X,\mu)$ has the standard meaning for $1\leq p\leq \infty$.  The space and measure will be suppressed where they are evident.  By a sequence $\{n_k\}$, we mean a sequence $\{n_k\}_{k=1}^{\infty}$ of natural numbers such that $n_k\to \infty$ as $k\to \infty$.  All spaces are endowed with a metric $d$ that gives rise to a locally compact, separable topology.  In particular, all spaces will be Hausdorff.

\section{Functional analytic properties of uniform rigidity}\label{s:mrah}
\begin{defin}
We say that a transformation $T$ of $X$ is {\bf uniformly rigid} if there exists a sequence $\{n_k\}$ such
that $d(T^{n_k}x,x)\to 0$ uniformly on $X$.  We say $T$ is {\bf rigid} if there exists a sequence $\{n_k\}$ such that for every finite measure $A\subset X$, $\mu(A\tR T^{n_k}A)\to 0$.
\end{defin}

We remark that the notions of rigidity and uniform rigidity often coincide.  See Lemma \ref{t:coinc}, for instance.  They are not the same in general.  For a trivial example, let $X$ be a set with countably many points given the discrete metric and a probability measure concentrated at one point.  Let $T$ be a transformation that is measure-preserving and is an infinite order shift on the remaining points.  This will be a transformation that is rigid but not uniformly rigid.  For a less trivial example, we will see that while it is possible to generically find weakly mixing, rigid, continuous transformations on a Cantor set, it is not possible for such a transformation to be uniformly rigid.  This is the content of Lemma \ref{t:mraah}.

\begin{lemma}\label{t:cont}
Assume $X$ is compact.
Let $f\in C(X)$, and define $f_k=f\circ T^{n_k}$.  If $T$ is uniformly rigid along $\{n_k\}$, then
$f_k\to f$ uniformly.  Conversely, if $f_k\to f$ uniformly for each $f\in C(X)$, then $T$ is
uniformly rigid.
\end{lemma}
Notice that in Lemma \ref{t:cont} we did not assume that $T$ is continuous.  The continuity of $f$ is essential in the proof, however.
\begin{proof}
Let $T$ be uniformly rigid and let $\eps>0$.  For each $x\in X$, choose $\delta$ such that $d(x,y)<
\delta$ implies that $|f(x)-f(y)|<\eps$.  Now, $f$ is continuous on a compact set and hence
uniformly continuous, so we may choose a $\delta>0$ that works for the whole space.  Choose such
a delta and a $k$ such that $d(T^{n_k}x,x)<\delta$ for all $x$.  Then, we have $|f_k(x)-f(x)|<\eps$.

For the converse, let $f^x(y)=d(x,y)$.  Then, $f^x(x)=0$ and $f_k^x(x)\to 0$, so that $T^{n_k}x\to
x$.  Now, for all $y\in X$, look at $B_{\eps/3}(y)$, and extract $y_1,y_2,\ldots,y_m$ that correspond
to a finite subcover.  For $x\in B_{\eps/3}(y_i)$, pick $k_i$ such that $|f^{y_i}_{k_i}(x)-f^{y_i}(x)|
<\eps/3$.  Then, since $d(x,y)<\eps/3$, we have that $d(T^{n_{k_i}}x,y)<2\eps/3$.  Now we see that
$d(T^{n_{k_i}}x,x)\leq d(x,y)+d(y,T^{n_{k_i}}x)<\eps$.  Choose $k=\max_i k_i$.
\end{proof}

We shall see later that when $T$ is ergodic, we cannot have this convergence uniformly in $f$, i.e. that it is impossible for $\{f_k\}_{k\in \bN}$ and all $f\in C(X)$ to be an equiconvergent class.

As $T$ is invertible, it induces a unitary representation of $\bZ$ on $L^2(X,\mu)$ by $U:f\mapsto f\circ T$.  The following result shows that uniform rigidity and rigidity have similar functional-analytic properties:

\begin{lemma}\label{t:conv}
Let $T$ be a uniformly rigid measure-preserving transformation on a compact space.
Then if $U$ is as above and $\{n_k\}$ is the uniform rigidity sequence, then $U^{n_k}\to I$ in the strong operator topology.
\end{lemma}

\begin{proof}
Let $C_c(X)$ denote the space of compactly supported continuous functions on $X$.  By the assumptions we have made about $X$ and $\mu$, it follows that $C_c(X)$ is dense in $L^2(X,\mu)$.  By Lemma \ref{t:cont}, $f_k\to f$ uniformly for each $f\in C_c(X)$.  If $f\in L^2$, choose $f^m\to f$ an $L^2$ approximation by continuous compactly supported functions.  We claim that $f_k\to f$ in $L^2$.  To show that $||f_k-f||_2$ is eventually smaller than $\eps$, it suffices to choose an $m$ such that $||f^m_k-f||_2<\eps/2$ and $||f^m_k-f_k||_2<\eps/2$.  The second approximation is obvious since $T$ is measure-preserving.  For the first one, note that we can assume that an arbitrarily large fraction of the measure is concentrated on a compact subset $Y$ of $X$ since $X$ is compact.  We can choose $m$ so that $\chi_Y\cdot f^m$ is within $\eps/4$ of $\chi_Y\cdot f$ in $L^2$ and then use the uniform convergence of $f^m_k\to f^m$ on $Y$ (by Lemma \ref{t:cont}) and the compact support of $f^m$ to get the desired approximation.
\end{proof}

By applying Lemma \ref{t:conv} to characteristic functions, it follows that uniform rigidity implies rigidity on compact spaces.  The fact that $U^{n_k} \rightarrow I$ in the strong topology applied to the characteristic function of $A$ says exactly $\mu (A \Delta T^{n_k}A)\rightarrow 0$.

It is easy to see that Lemma \ref{t:conv} does not have a converse as stated.  We also see that if $T$ is an arbitrary rigid measure-preserving transformation, then $U^{n_k}\to I$ in the strong operator topology.  However, uniform rigidity and rigidity do not coincide on compact spaces: for instance, the proof of Theorem \ref{t:mrah} exhibits a compact space where weak mixing and uniform rigidity are mutually exclusive, whereas there is a dense $G_{\delta}$ (in the appropriate space of representations) of continuous, weak mixing, and rigid transformations.

\begin{cor}\label{t:eigen}
Let $T$ be as before and let $\lambda$ be an eigenvalue of $U$.  Then $\lambda^{n_k}\to 1$.
\end{cor}

It seems that this is the strongest convergence we can expect from $U^{n_k}$:
\begin{prop}
Let $U$ be as before, with the assumption that $T$ is aperiodic.  Then, $U^{n_k}$ does not converge to the identity in the norm, nor do the Ces\`aro sums.
\end{prop}
\begin{proof}
It is a standard result from the spectral theory of dynamical systems (see \cite{MR1719722}, 3.4) that if $T$ is an aperiodic nonsingular transformation, then the spectrum of $U$ is dense in the unit circle, and is hence equal to the unit circle.  Let $\mathcal{A}$ be the commutative $C^*$-algebra generated by $U$ and its adjoint.  By the Spectral Theorem, $\mathcal{A}$ is isometrically isomorphic to the ring of complex continuous functions on the circle, and multiplication by $U$ is given by multiplication by $z$.  Suppose that $U^{n_k}\to I$ in norm.  Then, for any $\eps>0$ and $k$ sufficiently large, $||U^{n_k}-I||<\eps$.  It follows that the operator on $C(\bT)$ given by $f\mapsto z\cdot f$ converges to the identity along $\{n_k\}$, which is impossible even after passing to a subsequence.

As for the Ces\`aro sums, write
\[
f_k(z)=\frac{1}{k}\sum_{i=1}^k z^{n_i}.
\]
Viewing this expression as a complex-valued function on the circle $\bT$, we obtain
\[
f_k(\theta)=\frac{1}{k}\sum_{j=1}^k e^{in_j\theta}.
\]
Let $\chi_{S,k}$ be the characteristic function of the sequence $\{n_i\}$ as a subset of $\bN$, truncated to the $k^{th}$ term.  It is evident that
\[
\widehat{f_k}(n)=\frac{1}{k}\chi_{S,k}(n).
\]
As a function in $L^2(\bZ)$, $||\widehat{f_k}||_2=1/k$, so that the only candidate for the limit is $0$.  Since the Fourier transform is an isometric isomorphism between $L^2(\bT)$ and $\ell^2=L^2(\bZ)$, we have the claim.
\end{proof}

In general, it seems that it is difficult to expect norm convergence of operators arising in ergodic theory.  For some results along this line coming from uniform ergodic theory, see Corollary 2 in \cite{MR0417821}.  For instance, suppose that all mean ergodic operators on a Banach space are uniformly ergodic.  The mean ergodic operators are the ones that satisfy the existence of an operator $P$ such that \[Px:=\lim_{n\to\infty}\frac{1}{n}\sum_{k=1}^nT^kx\] exists for all $x$.  If the existence of $P$ implies that the C\`esaro sums of the $T^k$ converge to $P$ in norm, we say that $T$ is uniformly ergodic.  It is then true that the Banach space had to be finite-dimensional to begin with (see Corollary 3 in \cite{MR1867345}).

In contrast to the result of \cite{MR2032829} concerning rigidity, uniform rigidity does not appear to be generic in any measurable way.  It is a restrictive and nonmeasurable property in a precise sense:
\begin{lemma}\label{t:mraah}
Let $C$ denote Cantor space, equipped with a Borel probability measure. There exists no totally ergodic, measure-preserving system that is uniformly rigid in any metric compatible with the topology, except in the case where all the measure is concentrated at one point.
\end{lemma}

\begin{proof}
Any Cantor set is homeomorphic to a binary splitting tree, so we will show the result first in the case
of the standard metric on the $p$-adic ring $\bZ_2$ and then for any compatible metric.

For any fixed $r>0$, there are only $m<\infty$ many distinct balls of radius $r$.  By uniform rigidity and
the ultrametric inequality, it follows that we may choose a $k$ so that $T^{n_k}y\in B_r(x)$ if and
only if $y\in B_r(x)$.  Therefore, the balls of radius $r$ are each invariant under $T^{n_k}$,
contradicting total ergodicity unless the measure is concentrated in one such ball.  Iterating the
argument shows that then the measure must be concentrated at one point.

It remains only to show that for any compatible metric $d$ on $C$, that the proof carries over.  $C$
is metrizable and hence normal, so that if $C=\bigcup_{i=1}^m B_i$ is a union of balls,
we may separate $B_i$ and $B_j$
by neighborhoods for $i\neq j$.  For each $i$, $B_i$ is compact, so that $\inf\{d(x_i,x_j)\}>0$ for
$x_i\in B_i$ and $x_j\in B_j$.  By uniform rigidity and the ultrametric inequality,
we may choose $k$ so that for all $x\in X$, $d(T^{n_k}x,x)<\eps/2$, so that each of the $B_i$ is
$T^{n_k}$-invariant.
\end{proof}

Clearly the proof generalizes to any compact-open subset $C$ of $\bQ_p$.  Also, if $S$ is a transformation on $C$ and $(X,T)$ is another system such that $T\times S$ is totally ergodic on $X\times C$ in the product measure, then $T\times S$ is not uniformly rigid in the product metric.

We are now in the position to do the following:

\begin{proof}[Proof of Theorem \ref{t:mrah}]
The proof of the Theorem follows from the Jewett-Krieger Theorem (see \cite{MR833286}), which says that every finite measure-preserving ergodic transformation has a realization as a realization as a homeomorphism of a compact topological space.  Closer analysis of the proof reveals that the constructed space is actually a Cantor set.
\end{proof}

\section{Uniform rigidity and measurable weak mixing}
The following arises naturally from the results in \cite{MR1007412}.
\begin{quest}\label{t:conj}
Does there exists a finite measure-preserving ergodic system that is both weak mixing and
uniformly rigid?
\end{quest}
Notice that Glasner and Maon did not settle this question.  Their result is a reflection of the fact that the appropriate analogy of measure-theoretic rigidity in the topological context is uniform rigidity, because it has the right genericity property.  Furthermore, the notion of weak mixing they consider is topological weak mixing, which is again the correct analogous notion to measure-theoretic weak mixing.  The analogies are understood using representation-theoretic characterizations of both of these properties in measure-theoretic and topological contexts, where the analogy is clearer and easier to formulate.  See \cite{MR921351} for a the measurable representation-theoretic characterization.

One motivation for this question is that a (nontrivial) measure-preserving weakly mixing transformation that is uniformly rigid would yield an example of a measurably sensitive transformation that is not strongly measurably sensitive (see \cite{meas}.)  However, the existence of such transformations was shown by other methods in the same paper.  It was noted in the introduction that Glasner and Maon have shown that there exist many transformations that are uniformly rigid and whose only continuous eigenfunctions are constant.  We will be able to give some results in the direction of Question \ref{t:conj} in section 5.

\begin{lemma}\label{t:coinc}
Let $T$ be a continuous measure-preserving map on the unit interval $I$ or circle $S^1$ with a finite Radon measure $\mu$ that is nonatomic and gives positive measure to each nonempty open set.  Then rigidity and uniform rigidity are identical notions.
\end{lemma}
\begin{proof}
First, let $T$ be uniformly rigid, with uniform rigidity sequence $\{n_k\}$.  It suffices to show that for any open interval $I$ containing a point $x$ and $\delta>0$, there exists a $k$ such that $T^{n_k}(I)\subset B_{\delta}(I)$, where $k$ is uniform in $x$.  Notice that the image of an open interval will be an interval.  If we require $T$ to be uniformly rigid, then for some $k$ and all intervals $A$, $T^{n_k}(A)\subset
B_{\eps}(A)$.  It is clear that $A\tR T^{n_k}(A)$ can have at most two connected components, each
of which has diameter at most $\eps$.  A simple compactness argument shows that for each $r>0$ there
exists $\delta>0$ such that $\mu(B_{\delta}(x))<r$ for each $x$, whence the claim.

Conversely, suppose that $T$ is rigid along the same sequence.  A simple argument shows that for a fixed $\delta$, there exists $r>0$ such that $\mu(B_{\delta}(x))\geq r$ for every $x$.  Indeed, cut the interval into finitely many disjoint half-open balls (along with one closed ball at one end in the case of the interval) of size no more than $\delta/3$.  Let $r$ be the least measure of these balls, which is nonzero by assumption.  If $B$ is any ball of radius $\delta$, then $B$ contains one of the balls with measure at least $r$, so that we have the claim.  Seeing as $B\tR T^{n_k}(B)$ can have at most two connected components, we use rigidity to conclude that the diameter of these two components must go to zero.  In particular, we can cover $I$ with $\eps/2$-balls and pick a $k$ such that the diameter of each component of $B\tR T^{n_k}(B)$ is less than $\eps/2$ for every such ball $B$, so that for every point $x\in I$, $d(T^{n_k}x,x)<\eps$.
\end{proof}

There are no measure-preserving ergodic invertible continuous self-maps of the interval if all nonempty open sets have positive measure.  This can easily be seen as follows: we may assume $T(0)=0$ and $T(1)=1$ or $T(0)=1$ and $T(1)=0$.  Take a small closed interval $A$ containing $0$ with $1/2>m=\mu(A)>0$.  Then, $T$ maps $A$ to a closed interval containing $0$ or $1$.  Consider the full orbit of $A$ under all powers of $T$.  It is obvious that this a $T$-invariant subset of the interval.  Since open sets all have positive measure, the full orbit has measure no more than $2m<1$.  This follows from the observation that $T^2(A)=A$ by the continuity of $T$.

However, such maps do exist on $S^1$.  We do not see a way to extend the proof of Lemma \ref{t:coinc} to higher dimensional Euclidean spaces.  Indeed, the proof of the Lemma relies fundamentally on the notion that connected sets with large diameter have large measure, while no such relationship exists in dimensions $2$ and greater.  We can also say something about the relationship between rigidity and uniform rigidity on the interval even when the transformation is not continuous, obtaining a result that is reminiscent of Egorov's Theorem:

\begin{defin}
Let $(X,T)$ be a measurable transformation of a metric space, and let $Y\subset X$.  We say that $T$ is {\bf uniformly rigid on $Y$} if $d(y,T^{n_k}y)\to 0$ uniformly for $y\in Y$.
\end{defin}

\begin{thm}
Let $\mu$ be a finite Radon measure on a metric space $X$ that is nonatomic and gives positive measure to each nonempty open set.  Let $\eps>0$, and let $T$ be a rigid transformation on $X$.  Then there is a measurable subset $B\subset X$ such that $\mu(B)<\eps$ and $T$ is uniformly rigid on $X\setminus B$.
\end{thm}

\begin{proof}
Choose balls $\{B_{i}^{m}\}_{i=1}^{r_m}$ of radius at most
$2^{-m}$ so that
$$\mu\left(\left(\bigcup_{i=1}^{r_m}B_i^m\right)^c\right) <
2^{-m-1}\epsilon.$$
For each $m$, choose $n_{k_m}$ so that
$$\bigcup_{i=1}^{r_m}(B_i^m \Delta T^{n_{k_m}}B_i^m) < 2^{-m-1}\epsilon.$$
Then

$$B = \bigcup _{m=1}^{\infty}\left(\bigcup_{i=1}^{r_m}
\left(B_i^m\Delta T^{-n_{k_m}}B_i^m\right) \bigcup
\left(\left(\bigcup_{i=1}^{r_m}B_i^m\right)^c\right)\right)$$
has measure less than $\epsilon$ and $T$ satisfies
$$d(x, T^{n_{k_m}}x) < 2^{-m+1}$$ on $X\setminus B$.  Hence $T$ is uniformly rigid on $X\backslash B$.
\end{proof}

A natural question that arises from these results, aside from Question \ref{t:conj}, is the following:
\begin{quest}
In what situations do the notions of rigidity and uniform rigidity coincide?
\end{quest}

\section{Asymptotic convergence behavior}
It is well known that measurable weak mixing ($\mu(X)<\infty$) is mixing outside of some zero density collection of powers of the transformation (see the proof of Proposition \ref{t:d} for a precise statement.)  In the case of a weakly mixing transformation or an ergodic transformation with a sufficiently large set of eigenvalues, we will see that the phenomenon of uniform rigidity occurs with zero density.
It is well-known that if $(X,T)$ is a weakly mixing system, then $T$ is totally ergodic (see e.g. \cite{MR2371216}).

The following result is of independent interest, and might be helpful in resolving Question \ref{t:conj}.

\begin{prop}\label{t:d}
Let $(X,\mu)$ be a nonatomic measure space and let $T$ be a finite measure-preserving ergodic transformation.  Assume furthermore that $X$ is a compact metric space.  If $T$ is uniformly rigid, then the uniform rigidity sequence has zero density.
\end{prop}

\begin{proof}
Let $\{n_k\}$ be the uniform rigidity sequence.
By Lemma \ref{t:conv}, uniform rigidity sequences
are also rigidity sequences. In the finite measure-preserving case, the C\`esaro sums of $\mu(T^n(A) \cap A)$ converge to $\mu(A)^2$ and the limit of $\mu(T^{n_k}(A) \cap A)$ converges to $\mu(A)$.  As a
consequence, the density of $\{n_k\}$ is at most $\mu(A)$.  Since we can choose subsets of $X$ with arbitrarily small positive measure, this implies that the sequence $\{n_k\}$ has zero density.
\end{proof}

The proof of Proposition \ref{t:d} follows from rather soft principles which yield very little further insight into the nature of uniform rigidity sequences.  We thank the referee for posing the following natural question:
\begin{quest}
Which zero density sequences occur as uniform rigidity sequences for some system $(X,\mu)$ as above?
\end{quest}

\section{Group actions and generalized uniform rigidity}
Let $G$ be a countable group acting faithfully on a finite measure space $X$ by measure-preserving transformations.  By a faithful action, we mean that $1\neq g\in G$ implies $g$ does not act by the identity on $L^2(X)$.  Equivalently, the associated unitary action on $L^2(X)$ has trivial kernel.  Following Bergelson and Rosenblatt (see \cite{MR921351}, Theorem 1.9) the action of $G$ is weakly mixing if and only if the associated unitary representation of $G$ in $L^2(X)$ admits no nontrivial finite dimensional subrepresentations.  Thus, a $\bZ$-action is weakly mixing if and only if it admits no eigenvalues different from $1$, as $\bZ$ is abelian.  We remark for future reference that in the infinite measure-preserving case and the nonsingular finite measure case, we say that $T$ is weakly mixing if it admits no nonconstant $L^{\infty}$ eigenfunctions.

For a general nonabelian group, a finite dimensional subrepresentation is given by a function $f\in L^2(X)$ and finitely many group elements $\{g_1,\ldots,g_n\}$.  Let $V$ by the vector space on $\{g_i\cdot f\}$.  Then, the unitary representation of $G$, $\rho:G\to\mathcal{U}(L^2(X))$, admits a map $G\to GL(V)$ as a subrepresentation.

It follows that if $H<G$ and $Res^G_H\rho$ admits no finite dimensional subrepresentation, then the action of $G$ is weakly mixing.
\begin{defin}
Let $G$ be a countable group endowed with the discrete topology acting on a space $X$ as above.  The action of $G$ is {\bf uniformly rigid} if there exists a sequence $\{g_i\}$ of group elements that leaves every compact $K\subset G$, denoted $g_i\to\infty$, such that $d(x,g_i\cdot x)\to 0$ uniformly.
\end{defin}
We remark that though this is a natural condition for finitely generated groups, it makes sense for non-finitely generated groups.  Note that in this case for $g_i\to\infty$, it suffices for the $\{g_i\}$ to leave every finitely generated subgroup of $G$.  Most dynamical notions require the group to be locally compact and second countable, conditions that are automatically satisfied under our hypotheses.

\begin{lemma}
Let $G$ be a group acting on a space $X$ as above.  If $G$ contains a mixing automorphism, then the action of $G$ is weakly mixing.
\end{lemma}
\begin{proof}
See \cite{MR921351}, Corollary 1.6 or the discussion above, for example.
\end{proof}

For example, the usual $SL_2(\bZ)$-action on the torus is weakly mixing since every hyperbolic automorphism of the torus is mixing.  We can extend the standard $SL_2(\bZ)$-action on the torus to a (discontinuous but measure-preserving) $G$-action.

Let $A$ be any hyperbolic automorphism of the torus.  Now view the torus as $(\bT\times\bR)/\bZ$.  Let $X\subset \bT\times\bR$ be a fixed fundamental domain for the $\bZ$-action, which we choose to be a finite height cylinder.  Let $B_n$ be the map that glues the two ends of the cylinder together by a $1/n$-twist.  Explicitly, $B_n$ is given by $(x,y)\mapsto (x+y/n,y)\pmod 1$ on the domain $[0,1)\times [0,1)$.
Note that as $n$ gets large, the sequence of elements $\{B_n\}$ witnesses the uniform rigidity of the weakly mixing action of $\langle A,\{B_n\}\rangle$.  More generally:
\begin{thm}
Let $X$ admit a weakly mixing group action and a uniformly rigid action by nontrivial subgroups of a fixed group of automorphisms $G$.  Then there exists a $G$-action on $X$ that is simultaneously weakly mixing and uniformly rigid.
\end{thm}
\begin{proof}
We simply note that weak mixing is inherited from subgroups and uniform rigidity is inherited from infinite subgroups.  By finding continuous actions, one can obtain a continuous weakly mixing and uniformly rigid $G$-action.
\end{proof}

We remark briefly that there exist ergodic infinite measure-preserving uniformly rigid transformations as well as ergodic nonsingular finite measure uniformly rigid transformations.  Both can be seen by looking at a measure supported entirely on an orbit of an irrational rotation of the circle.  In the infinite measure-preserving case, we can take the measure to be counting measure.  It is obvious that this is ergodic and admits no $L^2$-eigenfunctions, but it admits $L^{\infty}$-eigenfunctions and is hence not weakly mixing.  In the nonsingular finite measure case, we choose a finite measure on the orbit such that the rotation is nonsingular.  Recall that if $T$ is a nonsingular transformation on a finite measure space $X$, then the action of $T$ on $L^2(X)$ is given by
\[
f\mapsto \sqrt{\omega}f\circ T,
\]
where $\omega$ is the Radon-Nikodym derivative of the measure with respect to $T$.  $T$ also acts on $L^{\infty}$ be precomposition.  In both of these cases, the action of $T$ gives rise to an isometry.
It is easy to see that $T$ is ergodic and admits no $L^2$-eigenfunctions but that there are nonconstant $\ell^{\infty}$-eigenfunctions, since $T$ yields an isometry of $\ell^{\infty}$ independently of the fact that $T$ is not measure-preserving.  Therefore, $T$ is not weakly mixing.  In both of the examples above, the associated unitary representations in $L^2$ admit no finite dimensional subrepresentations.

Having answered Question \ref{t:conj} in this context, we are left with the following:
\begin{quest}
To what extent is generalized uniform rigidity generic for ergodic group actions?
\end{quest}

\section{Acknowledgements}
The authors would like to thank the National Science Foundation for generous funding, as provided by REU Grant DMS-0353634, as well as Williams College for its funding and hospitable atmosphere during the 2006 SMALL REU program.  The second author would like to thank the 2008 SMALL REU program for its hospitality while the paper was put in final form.  The authors thank the referee for numerous comments that greatly improved the manuscript.


\begin{thebibliography}{99}
\bibitem[AS]{MR2032829}O. N. Ageev and C.E. Silva.  {\it Genericity of rigid and multiply recurrent infinite measure-preserving and nonsingular transformations}.  In {\it Proceedings of the 16th Summer Conference on General Topology and its Applications} (New York), volume 26, pages 357-365, 2001/02.
\bibitem[BR]{MR921351}Vitaly Bergelson and Joseph Rosenblatt.  Mixing actions of groups.  {\it Illinois J. Math.}, 32(1):65-80, 1988.
\bibitem[FLW]{MR1867345}Vladimir P. Fonf and Michael Lin and Przemys{\l}aw Wojtaszczyk. Ergodic characterizations of reflexivity of Banach spaces.  {\it J. Funct. Anal.} 187 (2003), no. 1, 146--162.
\bibitem[GM]{MR1007412}S. Glasner and D. Maon. Rigidity in topological dynamics.  {\it Ergodic Theory Dynam. Systems}, 9(2):309-320, 1989.
\bibitem[JKLSS]{meas}Jennifer James, Thomas Koberda, Kathryn Lindsey, Cesar E. Silva, and Peter Speh.  Measurable sensitivity.  {\it Proc. Amer. Math. Soc.}, 136(10):3549-3559, 2008.
\bibitem[L]{MR0417821}Michael Lin.  On the uniform ergodic theorem.  {\it Proc. Amer. Math. Soc.}, 43 (1974), 337--340.
\bibitem[N]{MR1719722}M. G. Nadkarni.  {\it Spectral theory of dynamical systems}.  Birkh\"auser Verlag, Basel, 1998.
\bibitem[P]{MR833286}Karl Petersen.  {\it Ergodic theory}, volume 2 of {\it Cambridge Studies in Advanced Mathematics}.  Cambridge University Press, Cambridge, 1983.
\bibitem[S]{MR2371216}C. E. Silva.  {\it Invitation to ergodic theory}, volume 42 of {\it Student Mathematical Library}.  American Mathematical Society, Providence, RI, 2008.
\end{thebibliography}
\end{document}